   \documentclass[12pt,final,onecolumn]{IEEEtran}

\IEEEoverridecommandlockouts              

\usepackage{amsfonts}
\usepackage{enumerate}
\usepackage{verbatim}
\usepackage{amssymb}
\usepackage{setspace}
\usepackage{epsfig}
\usepackage{url}
\usepackage{graphicx}

\usepackage{amsmath}

\usepackage[dvipsnames]{xcolor}

\def\urltilde{\kern -.15em\lower .7ex\hbox{\~{}}\kern .04em}
\def\urldot{\kern -.10em.\kern -.10em}

\def\urlhttp{http\kern -.10em\lower -.1ex\hbox{:}\kern -.12em\lower 0ex\hbox{/}\kern -.18em\lower 0ex\hbox{/}}

 \newtheorem{Theorem}{Theorem}
 \newtheorem{Lemma}{Lemma}

 \newtheorem{Assumption}{Assumption}
 \newtheorem{Problem}{Problem}
 \newtheorem{Corollary}{Corollary}
  \newtheorem{definition}{Definition}
 
 \newtheorem{Proposition}[Theorem]{Proposition}
 \newtheorem{remark} {Remark}
 
 \newtheorem{Example} {Example}

\newcommand*\diff{\mathop{}\!\mathrm{d}}

 \newcommand{\Int}{\operatorname{int}}

\newcommand{\diag}{\operatorname{diag}}

\newcommand {\R}{\mathbb R}

\newcommand{\N}{\mathbb{N}}

\newcommand{\be}{\begin{equation}}
\newcommand{\ee}{\end{equation}}
\usepackage{lineno}

\doublespace

 \title{\LARGE \bf   Entrainment to Subharmonic Trajectories  in \\
 Oscillatory  Discrete-Time~Systems\thanks{An abridged version of this paper has been submitted to the {\sl 27th Mediterranean Conference on Control and Automation (MED~2019)}. }}

\author{Rami  Katz, Michael Margaliot  and Emilia Fridman
\thanks{This research  was partially supported by  research grants from  the Israel Science Foundation    and the
 Binational Science Foundation.}
\thanks{RK and EF are  with School of Elec. Eng., Tel Aviv University, Israel.}
 \thanks{MM (Corresponding Author) is with the School of Elec. Eng. and the Sagol School of Neuroscience, Tel-Aviv
University, Tel-Aviv~69978, Israel.
{\tt\small E-mail: michaelm@eng.tau.ac.il}
}}

\begin{document}
\maketitle
\thispagestyle{empty}

\maketitle
\begin{abstract}
A matrix~$A$ is called totally positive~(TP) if all its minors are positive,
and  totally nonnegative~(TN) if all its minors are nonnegative. 
A square matrix~$A$ is called oscillatory if it is~TN and 
some power of~$A$ is~TP.
A linear time-varying system is called an
oscillatory discrete-time     system~(ODTS) if
 the matrix defining its evolution at each time~$k$
is oscillatory. We analyze the properties of 
$n$-dimensional time-varying nonlinear discrete-time systems 
whose variational system is an ODTS, and show that they 
have a well-ordered behavior.  More precisely, if the nonlinear system is
 time-varying and~$T$-periodic  then any trajectory either leaves any compact set or converges to an~$(n-1)T$-periodic trajectory, that is, a subharmonic trajectory. These results hold for any dimension~$n$. 
The analysis of such systems
requires establishing that 
a line  integral  of the Jacobian of the nonlinear system is an oscillatory matrix.
This is non-trivial, as the sum of two oscillatory  matrices is not necessarily oscillatory, and this carries over to integrals.  We derive several new sufficient
conditions guaranteeing
 that the line integral of a matrix is oscillatory, and 
demonstrate how this yields interesting classes of discrete-time
nonlinear systems that admit a well-ordered behavior. 
\end{abstract}


  \section{Introduction}

Positive dynamical systems arise naturally when the state-variables represent physical quantities
that can only take nonnegative 
values~\cite{farina2000,posi-tutorial}. 
For example, in compartmental systems the state-variables represent the ``density'' at each
 compartment~\cite{sandberg78},
in models of traffic flow or communication networks the state-variables represent the 
state of queues  in the system~\cite{shorten_tcp}, 
and in Markov chains the state-variables are probabilities~\cite{entrain_master_eqn}. 

Here, we introduce and analyze a new class  of 
 positive systems called \emph{oscillatory discrete-time systems}. 
Recall that 
a matrix~$A\in\R^{n\times m}$ is called \emph{totally positive}~(TP)
if every minor of~$A$ is positive,
and  \emph{totally nonnegative}~(TN) if every minor of~$A$  is non-negative.\footnote{Unfortunately, the terminology in this field is not uniform. We follow the terminology used  in~\cite{total_book}.} 
TN and~TP matrices have a remarkable variety
of interesting mathematical properties~\cite{total_book,pinkus}.
One important property is that multiplying a vector by a~TP
 matrix cannot increase the number of sign variations in the vector. 

Oscillatory matrices are in the ``middle ground'' between TN and~TP matrices.
A matrix~$A\in\R^{n\times n}$ is called \emph{oscillatory} 
if~$A$ is TN and there exists an integer~$k>0$ such that~$A^k$ is~TP. 
For example, it is easy to verify that all the minors of
\be\label{eq:aosc}
A:=\begin{bmatrix}  0.2& 0.1 & 0\\9&11&1\\0 &1 
&3    \end{bmatrix} 
\ee
 are nonnegative so~$A$ is (TN) (but not~TP as it has zero entries), and also
 that all the minors of~$A^2=\begin{bmatrix}  0.94& 1.12 & 0.1
\\100.8&122.9&14\\9 &14 &10    \end{bmatrix}$ are positive, so~$A$ is oscillatory.

The product of two TP/TN/oscillatory matrices is a  TP/TN/oscillatory matrix, but the sum
 of two TP/TN/oscillatory  matrices
is not necessarily a TP/TN/oscillatory matrix. 
For example,   the matrix~$A=\begin{bmatrix}  
1 & 0.1 \\ 9 & 1
\end{bmatrix}$
and its transpose~$A'$  are~TP (and thus in particular TN and oscillatory), yet
$
A+A'=\begin{bmatrix} 
2 & 9.1 \\ 9.1 & 2
\end{bmatrix}
$
is not~TN (and thus not~TP nor oscillatory), as~$\det(A+A')<0$.

TP matrices have important applications in the asymptotic  analysis of both continuous-time
and discrete-time dynamical systems. 
Schwarz~\cite{schwarz1970} introduced the notion of a \emph{totally positive 
differential system~(TPDS)}. This is the   
 linear time-varying~(LTV) system
\be\label{eq:lopo}
\dot x(t)=A(t)x(t),
\ee
 satisfying  that 
the associated \emph{transition matrix}~$\Phi(t_1,t_0)$ is~TP for any pair~$(t_1,t_0)$ 
with~$t_1>t_0$. 
The transition matrix is the matrix satisfying~$x(t_1)=\Phi(t_1,t_0)x(t_0)$ for 
all~$x(t_0)\in\R^n$. In the particular case where~$A(t)\equiv A$ the transition
 matrix is~$\Phi(t,t_0)=\exp((t-t_0)A)$, and then~\eqref{eq:lopo} is TPDS iff~$A$ is tridiagonal
with positive entries on the super- and sub-diagonals. 
Schwarz used the~VDP to show that the number of 
sign variations in~$x(t)$
is a (integer-valued) 
Lyapunov 
 function for the TPDS~\eqref{eq:lopo}.
It was recently shown~\cite{fulltppaper} that
TPDSs  have 
 important applications in the stability analysis of continuous-time nonlinear cooperative dynamical systems with a tridiagonal Jacobian.

An extension to   discrete-time systems, called 
a \emph{totally positive  discrete-time  system}~(TPDTS),
 has been suggested recently~\cite{rola}. 
The  LTV
\be\label{eq:ltvds}
x(k+1)=A(k)x(k), 
\ee
with~$A:\N\to \R^{n\times n }$,
 is called a~TPDTS  
if~$A(k)$ 
 is~TP for all~$k\in\mathbb{N}$. 
It was shown that time-varying
nonlinear systems, whose variational equation is a~TPDTS,
 satisfy strong asymptotic properties  including entrainment to a periodic excitation. 
The variational equation is an~LTV with a matrix
described by a  line integral of the Jacobian of the nonlinear system.  
Since the sum of two~TP matrices is not necessarily~TP, it is
not trivial to verify that this  line integral   is indeed~TP.  

The main contributions of this paper are two-fold.
First, we introduce the new notion of an \emph{oscillatory discrete-time system}~(ODTS).
The  LTV~\eqref{eq:ltvds}
  is called an~ODTS  
if~$A(k)$ 
 is oscillatory  for all time~$k$. This is an important generalization of a~TPDTS,
as oscillatory matrices are much more common than~TP matrices. 
We analyze the properties of  discrete-time time-varying
nonlinear systems, whose variational equation is 
an~ODTS, and show that they
 satisfy useful asymptotic properties. In particular, if the~$n$-dimensional
 time-varying nonlinear system is~$T$-periodic then every solution either leaves every compact set or converges to an~$(n-1)T$-periodic solution, i.e. a subharmonic solution. 

The variational equation associated 
with the nonlinear system 
is an~LTV with a matrix
described by a  line integral of the Jacobian of the nonlinear system.  
Since the sum of two oscillatory  matrices is not necessarily oscillatory, it is
not trivial to verify that this  line integral   is indeed oscillatory.  

The  second contribution of this paper is deriving 
several new sufficient conditions guaranteeing that the line 
 integral   of a matrix is oscillatory. Our first condition considers
 the  special case of a system with  scalar nonlinearities.
 In this case we show that 
   the integration can be performed in closed-form. 
The other conditions  are based on   sufficient conditions for a matrix to be oscillatory or~TP. 
We demonstrate how these conditions yield new classes of discrete-time 
nonlinear systems 
with a well-ordered behavior.

The remainder of this paper is organized as follows: Section~\ref{sec:Prelim} reviews known definitions and results that will be used later on including the~VDPs of~TN and~TP matrices,
and~TPDTSs. 
The next two sections describe our main results. Section~\ref{sec:odts}
defines and analyzes~ODTSs. Section~\ref{sec:scon}
 provides several sufficient  conditions verifying that the line integral of the Jacobian of a time-varying
nonlinear system is oscillatory. This section also details several applications of the theoretical results.  
The final section concludes and describes several topics for further research. 

We use standard notation. The set of nonnegative integers is~$\N:=\{0,1,2,\dots\}$. 
Matrices [vectors] are denoted by capital [small] letters.
The transpose of a matrix~$A$ is denoted~$A'$.
We use~$\diag(v_1,\dots,v_n)$ to
denote the~$n\times n$  diagonal  matrix with   
entries~$v_1,\dots,v_n $ on the diagonal.

\section{Preliminaries}\label{sec:Prelim}
We begin by reviewing the VDP of TN and~TP matrices. More details and proofs can be found in the excellent monographs~\cite{total_book,pinkus,gk_book}.
For a vector~$z\in\R^n$ with no zero entries the number of sign variations
 in~$z$ is
\be\label{eq:sigmdfr}
\sigma(z):=\left|\{i \in \{1,\dots,n-1\} : z_i z_{i+1}<0\} \right | .
\ee
For example, for~$n=3$ consider  the vector~$z(\varepsilon):=\begin{bmatrix}  2 & \varepsilon & -3 \end{bmatrix}' $.
 Then for \emph{any}~$\varepsilon  \in \R\setminus\{0\}$,  $\sigma(z(\varepsilon))$   is well-defined and equal to
one. 
More generally, the domain of definition of~$\sigma$    can be extended, via continuity, to the  set:
\begin{align*}
V := & \{z\in\R^n:   z_1 \not =0,\; z_n \not=0,\; \text{and if }  z_i=0 
\\
&\text{ for some~$i \in \{2,\dots,n-1\} $ then } z_{i-1} z_{i+1}<0\}.
\end{align*}
  
	We   recall   two more   definitions for   
		the number of sign variations in a 
		vector~\cite{total_book} that are well-defined for any~$y\in\R^n$. 
 	Let \[ s^-(y):=\sigma(\bar y), \] 
	where~$\bar y$ is the vector obtained from~$y$ by deleting all zero entries,
	and let \[s^+(y):=\max_{x\in P(y)} \sigma(x) , \]
	where~$P(y)$ is the set of all vectors obtained by replacing every zero entry of~$y$ by either~$-1$ or~$+1$. 
	For example,
	for~$y=\begin{bmatrix} -1& 0 &0 &4  \end{bmatrix}'$, 
	$s^-(y)=1$ and~$s^+(y)=3$. 
These definitions imply that 
	\be\label{eq:implthat}
 0\leq s^-(y) \leq s^+(y) \leq n-1  \text{ for all } y\in\R^n.
\ee
An   important observation 
is that~$s^-(y)=s^+(y)$ iff~$y\in V$.

	A classical result~\cite{total_book}   states
	that if~$A\in \R^{n\times m}$ is~TP 
	then
	\[
	s^+(Ax)\leq s^-(x) \text{ for all } x \in \R^m \setminus \{0\},
	\]
	whereas if~$A$ is~TN (and in particular if it is~TP) 
					then
	\be\label{eq:vdptn}
					s^-(Ax)\leq s^{-}(x) \text{ for all }x\in \R^m. 
	\ee
These are the VDPs of TP and TN matrices. For example,
the matrix~$A=\begin{bmatrix} 1 &2 \\ 1&4\end{bmatrix}$ is~TP
and for~$x=\begin{bmatrix} 1 &-1 \end{bmatrix}'$, we have
\[
s^+(Ax)=s^+(\begin{bmatrix}-1&-3\end{bmatrix}') < s^-(x).
\]

For  square  
matrices (which is the relevant case when considering the
transition matrices of 
dynamical systems) more precise results are known. 
Recall that a matrix is called \emph{strictly sign-regular of order~$k$} 
(denoted~$SSR_k$) if
   its minors of order~$k$ are either all positive or all negative. 
For example,~$A=\begin{bmatrix}  1& 2 \\ 3&4 \end{bmatrix} $ 
is~$SSR_1$ because all its entries are positive, and~$SSR_2$
because it single minor of order~$2$ is negative. 
It was recently shown~\cite{CTPDS} that if~$A\in\R^{n\times n}$ 
is non-singular
then for any~$k\in \{1,\dots,n-1\}$ 
we have that~$A$ is~$SSR_k$  iff
\[
			x \in\R^n\setminus\{0\} \text{ and }	s^-(x) \leq k-1 \implies 	s^+(Ax)\leq k-1.
\]

For example, for~$k=1$ this implies that for a non-singular matrix~$A\in\R^{n\times n}$ the following
   two conditions are equivalent: (1)~all the entries of~$A$ are either all positive or all negative;
and (2)~for every~$x\not =0$ with all entries non-negative or all
non-positive the vector~$Ax$ has all entries positive or all negative. 

We now review applications of total positivity to discrete-time
dynamical systems.

\subsection{Totally Positive Discrete-Time  Systems} 
Consider the  discrete-time  LTV~\eqref{eq:ltvds}
with~$A:\N\to \R^{n\times n }$. 
The system 
 is called a~TPDTS~\cite{rola}
if~$A(k)$ is~TP for all~$k\in \N$.  
Intuitively speaking, this 
 is the discrete-time
analogue  of a~TPDS.
The~VDP and~\eqref{eq:implthat}
imply that for any~$x(0)\in\R^n\setminus\{0\}$  we have 
\be\label{eq:dec}
  \dots     \leq s^-(x(1))\leq  s^+(x(1)) \leq s^-(x(0))\leq s^+  (x(0)).
\ee
In other words, both~$s^-(x(k))$ and~$s^+(x(k))$ can be viewed as
 integer-valued Lyapunov functions for the trajectories  of a~TPDTS.  
Furthermore, there can be no more then~$n-1$ strict inequalities in~\eqref{eq:dec},
 as~$s^-$ and~$s^+$ take values in~$\{0,1,\dots,n-1\}$. 
This implies that there exists~$m\in \N$ such that~$s^-(x(k))=s^+(x(k))$
 for all~$k\geq m$, i.e.~$x(k) \in V$ for all~$k\geq m$. In particular,~$x_1(k) \not =0$ (and~$x_n(k) \not =0$) for all~$k\geq m$. 
Moreover, it was  shown in~\cite{rola} that   there exists~$p\in\N$  such that 
 the following \emph{eventual monotonicity}
property holds: 
  either~$x_1(k )>0$ for all~$k\geq p$ 
	or~$x_1(k )<0$ for all~$k\geq p$ 
	(and similarly for~$x_n(k)$).

This property  can be applied to study the asymptotic properties
 of time-varying
\emph{nonlinear} discrete-time systems. 
Consider the   system
\be\label{eq:fsys}
x(k+1)=f(k,x(k)).
\ee
We assume that~$f:\N\times \R^n\to\R^n$ is~$C^1$ with respect
 to its second variable, and
 denote its Jacobian  by~$J(k,x):=\frac{\partial}{\partial x} f(k,x)$. 
We also assume that the trajectories of~\eqref{eq:fsys}
evolve on a compact and convex state-space~$\Omega\subset \R^n$.
For~$a \in\Omega$ and~$j\in\N$, let~$x(j,a)$ denote the solution of~\eqref{eq:fsys} 
at time~$j$ with~$x(0)=a$.

Fix~$a,b \in \Omega$ and let~$z(k):=x(k,b)-x(k,a)$.
Then (see, e.g.~\cite{rola})
\be\label{eq:zsys}
z(k+1)=M(k) z(k),
\ee
where
\begin{align}\label{eq:yassa}
M(k)=M(k,a,b):=\int_0^1 J(k,rx(k,b)+(1-r)x(k,a) ) \diff r.
\end{align}
The LTV system~\eqref{eq:zsys} is called 
the \emph{variational equation} associated with~\eqref{eq:fsys}, as
it describes how the variation between the 
two solutions~$x(k,b)$ and~$x(k,a)$
evolves in time. 

We pose two assumptions.
\begin{Assumption}\label{assume:jaco}
The matrix
\be\label{eq:deffmat}
F(k,a,b):=\int_0^1 J(k,ra+(1-r)b ) \diff r
\ee
is~TP for all~$k\in\N$ and all~$a,b\in \Omega$.
\end{Assumption}
Note that this implies that~\eqref{eq:zsys} is a~TPDTS. 

\begin{Assumption}\label{assume:tper}
There exists~$T \in \{1,2,\dots\}$ such that 
the map  in~\eqref{eq:fsys}
 is~$T$-periodic, that is, 
\[
f(k,a)=f(k+T,a) \text{ for all }k\in \N \text { and all } a \in \Omega. 
\]
 \end{Assumption}
Note that in the particular case where~$f$ is time-invariant   this holds
(vacuously) for every~$T \in \mathbb{N}$.

\begin{Theorem}\cite{rola} \label{thm:mrola}
If  Assumptions~\ref{assume:jaco}
  and~\ref{assume:tper} hold then 
	 every solution of~\eqref{eq:fsys} emanating from~$\Omega$ 
converges to a $T$-periodic solution of~\eqref{eq:fsys}.
\end{Theorem} 

 If the time-dependence in~$f$ is due to an input (or excitation)~$u$, that is, $f(k,x(k))=g(u(k),x(k))$ for some map~$g$ then Assumption~\ref{assume:tper}
holds if~$u$ is $T$-periodic.   Thm.~\ref{thm:mrola} then implies
 that the system \emph{entrains} to the periodic excitation, as every solution converges to a periodic solution with the same period~$T$. Entrainment is an important property in many natural and artificial systems~\cite{2017arXiv171007321M,RFM_entrain,entrain2011}.
For example, many biological processes in our bodies, like the sleep-wake cycle,  entrain to the 24h-periodic   solar day.

In the special case where~$f$ is time-invariant Thm.~\ref{thm:mrola}
yields the following result. 
\begin{Corollary}\label{coro:tinv}\cite{rola} 
Consider the time-invariant nonlinear system
\be\label{eq:fsysti}
x(k+1)=f( x(k)) 
\ee
whose trajectories evolve on a compact and convex state-space~$\Omega \subset \R^n$. 
Suppose that  
\be\label{eq:defab}
F(a,b):=\int_0^1 J( ra+(1-r)b ) \diff r
\ee
is~TP for   all~$a,b\in \Omega$. Then every solution of~\eqref{eq:fsysti} emanating from~$\Omega$
converges to an equilibrium point. 
\end{Corollary} 

Note that the equilibrium point is not necessarily unique.

The condition on~$F(a,b)$ implies
  that every minor of~$J(x)$ is positive
 for all~$x\in \Omega$. In particular, the first-order minors, i.e. the entries of~$J(x)$
are positive, so  the nonlinear system is strongly 
cooperative~\cite{smith_survey_2017,hlsmith}. 
 The conditions here require more than strong cooperativity
 and as a consequence yield more powerful results on the asymptotic behavior of the
 system (see, e.g.~\cite{smith_planar,coop_mini_review}).

In the particular    case of planar systems, the conditions here
require that the entries of~$J(x)$ are positive,
and that~$\det J(x)$ is positive. 
The latter condition is known to be an orientation-preserving condition
that has been used in the analysis of planar cooperative systems~\cite{smith_planar}.

The next result, which seems to be new, 
 shows that total positivity (in fact, a slightly weaker condition)
 implies an orientation-preserving property (with respect to 
 a specific  
 order)
 for \emph{any} dimension~$n$.
For two vectors~$x,y\in\R^n$,   we write~$x\ll y$ if~$x_i<y_i$ 
for all~$i\in\{1,\dots,n\}$. 
Let~$D_\pm\in\R^{n\times n}$ be the diagonal matrix with~$d_{ii}=(-1)^{i+1}$ for
all~$i\in\{1,\dots,n\}$. 
Note that~$(D_\pm)^{-1}=D_\pm$. We say that~$z\in\R^n$ is \emph{alternating}
if~$z_i z_{i+1}<0$ for all~$i\in\{1,\dots,n-1\}$. 
This implies of course that~$s^-(z)=s^+(z)=n-1$. 


\begin{Lemma}\label{lem:nrety}
Let~$P\in\R^{n\times n}$ be~TN and nonsingular.
 If~$x,y\in\R^n$ are such that
\be\label{eq:pozie}
D_\pm P D_\pm x \ll  D_\pm P D_\pm y
\ee 
then
\[
  x \ll   y . 
\]
\end{Lemma}
The proof is placed in the Appendix. 

\begin{Example} 
Consider the~TP matrix~$P=\begin{bmatrix} 1&2\\ 3&8\end{bmatrix}$.
Then~\eqref{eq:pozie} becomes
$
					\begin{bmatrix} 1&-2\\ -3&8\end{bmatrix} (x-y)\ll 0 	
$ 
and this holds iff~$x_1 - y_1<0$ and~$\frac{3}{8}<\frac{x_2-y_2}{x_1-y_1}<\frac{1}{2}$, so in particular~$x\ll y$. 
\end{Example}

In the context of the LTV~$z(k+1)=P z(k)$, $z(0)=z_0\in\R^n$, 
Lemma~\ref{lem:nrety} implies the following. Suppose that~$P$ is~TN and non-singular  and let~$y(k):=D_\pm z(k)$. 
Then    it is \emph{not} possible that for some~$i\geq 1$ we have
\be\label{eq:mocpyte}
y(0)\ll y(1)\ll \dots\ll y(i-1) \ll y(i ) \text{ and } y(i )  \gg y(i+1).
 \ee
Indeed, the last inequity here yields
\[
     D_\pm P D_\pm    y(i-1) \gg  D_\pm P D_\pm  y (i),
\]
so Lemma~\ref{lem:nrety}  gives
\[
         y(i-1) \gg     y (i),
\]
and this contradicts~\eqref{eq:mocpyte}.

     Smillie~\cite{smillie} and Smith~\cite{periodic_tridi_smith} proved 
convergence to an equilibrium and entrainment 
in a  certain class  of continuous-time 
nonlinear dynamical systems. Their results are based on using  
the number of sign variations in the solution of the associated (continuous-time)
variational system
 as an integer-valued Lyapunov function. 
It was recently shown that these results are closely
 related to the theory of TPDSs~\cite{fulltppaper}. 
Thm.~\ref{thm:mrola} and Corollary~\ref{coro:tinv} may be regraded as 
discrete-time analogues of these results.

It is well-known that asymptotically stable linear systems entrain  to periodic excitations. However,  nonlinear systems do not necessarily entrain. This is true even for strongly monotone systems.  
Ref.~\cite{Takac1992} provides interesting 
examples of \emph{continuous-time}, 
strongly  cooperative dynamical systems 
whose vector field is~$T$ periodic 
and admit a
solution  that is  periodic with \emph{minimal} period~$nT$, for any integer~$n\geq 2$. 
Furthermore, this subharmonic solution may be asymptotically stable.

In order to apply Thm.~\ref{thm:mrola} and Corollary~\ref{coro:tinv} one needs to verify that 
  the line integral of the Jacobian is~TP. 
This is not trivial because the sum of two~TP matrices is not necessarily a~TP matrix,
and this    is naturally carried over to integrals. 
\begin{Example}\label{exa:nirt}
It is straightforward to verify that~$A(t) =\begin{bmatrix}
1.01 & t+1\\\frac{1}{t+1}& 1 
\end{bmatrix} $ 
 is~TP for all~$t\in[0,1]$, yet 
$
			\int_0^1 A(t) \diff t=  \begin{bmatrix}
1.01 & 3/2 \\ \ln(2) & 1 
\end{bmatrix}$ 
  is 
not~TP (and not even~TN), as it has a negative determinant. 
\end{Example} 

A matrix $A\in\R^{n\times n}$ is called \emph{oscillatory} if it is~TN and there exists~$k\in \mathbb{N}$ such 
that $A^k$ is~TP. The smallest such~$k$ is called the \emph{exponent} of the oscillatory matrix~$A$.
Oscillatory matrices are in the ``middle ground''
between~TN and~TP matrices, and
are much more common than~TP matrices in applications.
Indeed, it is well-known that a~TN matrix~$A\in\R^{n\times n}$ is oscillatory if and only if it 
		is non-singular and irreducible~\cite[Ch.~2]{total_book}, and that in this case~$A^{n-1}$ is~TP.  
		The next example demonstrates this. 
\begin{Example}\label{exa:osc}
Consider the    tridiagonal   matrix
		\be\label{eq:trida}
							A=\begin{bmatrix} 
																	a_1 & b_1 & 0 & \dots & 0											 \\
																	c_1 & a_2 & \ddots & \dots & \vdots											 \\
																	0 & \ddots & \ddots & \dots & \vdots											 \\
																	\vdots & \ddots & \ddots & \dots & b_{n-1}											 \\
																	0 & \dots & \dots & c_{n-1} & a_{n}											 
									\end{bmatrix} 
		\ee
		with~$ b_i,c_i\geq0$  for all~$i$.
		In this case, the \emph{dominance condition}
		\be\label{eq:domi}
		a_i \geq b_i+c_{i-1} \quad \text{for all } i\in\{1,\dots,n\},
		\ee
		with~$c_0:=0$ and~$b_n:=0$, guarantees that~$A$ is TN (see e.g.~\cite[Ch.~0]{total_book}). If, furthermore,~$b_i,c_i>0$ for all~$i$ then~$A$ is irreducible. Thus, if~$A$ is also non-singular then it is oscillatory. 
		\end{Example}

The next two sections describe our main results. 

\section{Oscillatory Discrete-Time Systems} \label{sec:odts}
We begin by introducing the new notion of an~ODTS. 
\begin{definition}
The discrete-time LTV 
\begin{equation}\label{eq:OscSys}
y(k+1)=A(k)y(k),  
\end{equation}
with~$A:\mathbb{N}\to\R^{n\times n}$, is called an \emph{ODTS of order~$p$}
 if $A(k)$ is~oscillatory for all~$k\in \mathbb{N}$, and every product of~$p$ matrices in the form:
\[
											A(k_p) \dots A(k_2)  A(k_1),      \quad 0\leq k_1< \dots <k_p,
\]
is TP.
\end{definition}

For example, if~$A(k)$ is~TP for all~$k$ then~\eqref{eq:OscSys} is an~ODTS of order~$1$.
Also, since the product of any~$n-1$ oscillatory matrices
 is~TP~\cite{pinkus},
 \eqref{eq:OscSys} is always an ODTS of order~$n-1$.

We now describe the  applications of ODTS to the 
 time-varying nonlinear system:
\begin{equation}\label{eq:nonlinOsc}
x(k+1) = f(k,x(k)),
\end{equation}
where $f(k,x)$ satisfies Assumption~\ref{assume:tper}. We assume that the trajectories of \eqref{eq:nonlinOsc} evolve in a compact and convex state-space~$\Omega\in \mathbb{R}^n$. 
For~$k\in\N$ and~$a,b\in \Omega$,
let
\[
F(k,a,b):=\int_0^1 J(k,ra+(1-r)b ) \diff r.
\]

We pose the following assumption.
\begin{Assumption}\label{assume:jaco1}
The   system
\be\label{eq:deffmat1}
z(k+1)=F(k,\cdot,\cdot)z(k) 
\ee
is an ODTS of order~$h $.
\end{Assumption}

We can now state the main result in this section. 
\begin{Theorem}\label{Thm:Oscillsys}
Suppose that Assumptions~\ref{assume:tper}  and~\ref{assume:jaco1} hold.
Let~$u:=h T$.  
 Then every solution of \eqref{eq:nonlinOsc} emanating from~$\Omega$ 
converges to a $u$-periodic solution of \eqref{eq:nonlinOsc}.
\end{Theorem}

\begin{remark}
If~$F(k,a,b)$ is~TP for all~$k\in\N$ and all~$a,b\in \Omega$
then Assumption~\ref{assume:jaco1} holds with~$h=1$ so
Thm.~\ref{Thm:Oscillsys} implies that every solution of~\eqref{eq:nonlinOsc} emanating from~$\Omega$ 
converges to a $T$-periodic solution of \eqref{eq:nonlinOsc}. This recovers 
  the~TPDS case.
If~$F(k,a,b)$ is oscillatory   for all~$k\in\N$ and all~$a,b\in \Omega$
then in particular every product of~$n-1$ matrices is~TP, 
so Thm.~\ref{Thm:Oscillsys} implies that every solution of~\eqref{eq:nonlinOsc} emanating from~$\Omega$ 
converges to an~$(n-1)T$-periodic solution of \eqref{eq:nonlinOsc}.
\end{remark}

\begin{remark}
The LTV~\eqref{eq:OscSys} is of course a special case of~\eqref{eq:nonlinOsc} with Jacobian~$J(k,x(k))=A(k)$, and 
thus~$F(k,a,b)=A(k)$ for all~$a,b\in\Omega$ and all~$k\in\mathbb{N}$. 
We conclude that if~$A(k)=A(k+T)$ for all~$k\in\mathbb{N}$
then every solution of an~ODTS  of order~$h$
converges to periodic solution of~\eqref{eq:OscSys}
with period~$u:=h T$.   
\end{remark}

\begin{IEEEproof}[Proof of Thm.~\ref{Thm:Oscillsys}]
Pick~$\alpha,\beta \in \Omega$ with~$\alpha\neq \beta$.
Let
\[
z(k):=x(k,\beta)-x(k,\alpha) 
\]
 and recall that~$z$
satisfies the  variational equation~\eqref{eq:zsys}, with
\begin{equation*}
M(k,\alpha,\beta):=\int_0^1 J(k,rx(k,\beta)+(1-r)x(k,\alpha) ) \diff r.
\end{equation*}
Assumption~\ref{assume:jaco1} implies that~$M(k,\alpha,\beta)$
is oscillatory.
Let~$v(k):=z( ku )$. Then
\begin{equation}\label{eq:subseq}
v(k+1) = M((k+1)u-1,\alpha,\beta) \dots  M(k u ,\alpha,\beta) v(k).
\end{equation}
The  product on the right-hand side
 includes~$u=hT$ matrices, and is~TP, as 
  the product of any~$h$ matrices is~TP, and   the product of 
any two~TP matrices is~TP.
 Thus,~\eqref{eq:subseq}
is a~TPDS. Thm.~6 in~\cite{rola}  implies
the following  eventual monotonicity property: 
there exists $m\in \mathbb{N}$ such that either $v_1(k)>0$ for all $k\ge m$  or $v_1(k)<0$ for all $k\ge m$.

Pick~$a\in \Omega$ and let $x(k,a)$ denote the trajectory of \eqref{eq:nonlinOsc} emanating from~$a$. 
Let~$b:=x(u,a)$.
If~$x(k,a)$ is~$u$-periodic then
 there is nothing to prove. Therefore, we can assume that the trajectories $x(k,a)$ and~$x(k,b)$
are not identical. Note that Assumption~\ref{assume:tper} implies that both trajectories are solutions of~\eqref{eq:nonlinOsc}. By the eventual monotonicity property,
  there exists~$m\in \mathbb{N}$ such that, without loss of generality,
\begin{eqnarray}\label{eq:EvMono}
x_1((k+1)u,a)-x_1(ku,a)>0 \text{ for all } k\geq m.
\end{eqnarray}

Let
\begin{align*}
\omega_{u}(a) :=\{& p: \text{ there exist } m_i\in\mathbb{N}  \text{ with }
  m_1<m_2<\dots \\& \text{ such that }\lim_{i\to\infty}x(m_i u,a)=p  \},
\end{align*}
that is,   the~$u$-omega limit set corresponding to~$a$.
 By compactness of $\Omega$ it follows that $\omega_{u}(a)\neq \emptyset$. We now show 
that~$\omega_{u}(a)$  is a singleton. Assume that there exist~$p,q\in \omega_{u}(a)$,
with~$p\neq q$.
 We claim that $p_1=q_1$. Indeed, there exist sequences~$\left\{m_k\right\}_{k=1}^{\infty}$ and~$\left\{s_k\right\}_{k=1}^{\infty}$ such that~$p = \lim_{m_k \to \infty}x(m_k u,a)$ and~$q =\lim_{s_k \to \infty}x(s_ku,a)$. Passing to sub-sequences, if needed, we
 may  assume that~$m_k<s_k<m_{k+1}$ for all~$k\in \mathbb{N}$. Now~\eqref{eq:EvMono} implies that~$p_1=q_1$. We conclude that any two points in~$\omega_{u}(a)$   have the same first coordinate. 

Consider the trajectories emanating from~$p$ and~$q$, that is,~$x(k,p)$ and $x(k,q)$. Since $p,q\in \omega_{u}(a)$ and $\omega_{u}(a)$ is an invariant set, \[
x(ku,p),x(ku,q)\in \omega_{u}(a)  \text{ for all } k\in \mathbb{N}.
\]
 This implies that 
\begin{equation*}
x_1(ku,p)=x_1(ku,q) \text{ for all }  k\in \mathbb{N}.
\end{equation*}
However, this contradicts the eventual monotonicity of~\eqref{eq:subseq}.
We conclude that~$\omega_u(a)$ is a singleton, and this completes the proof.  
\end{IEEEproof}

\begin{remark}
Note that the proof of Thm.~\ref{Thm:Oscillsys} relies on the fact that any  product of~$u=hT$ matrices in~\eqref{eq:subseq} is~TP. In practice,  it may be the case that a product of a smaller
number of    matrices in the variational equation is~TP. In this case, every
 solution~$x(k,a)$ will converge to a periodic 
solution of~\eqref{eq:nonlinOsc} with period less than~$hT$. Nevertheless, the minimal period of the limit solution must divide~$hT$.
\end{remark}

 The next subsection provides several sufficient conditions guaranteeing
that Assumption~\ref{assume:jaco1}  indeed holds, and applications to several dynamical systems.

\section{Conditions guaranteeing that a matrix 
line integral is oscillatory}\label{sec:scon}

Our first sufficient condition is based on the sufficient
condition for a matrix to be oscillatory 
 described in Example~\ref{exa:osc}.

\subsection{Discretizing nonlinear tridiagonal strongly cooperative~systems}
Consider the nonlinear time-varying dynamical system~$\dot x=f(t,x)$.
Let
\be\label{eq:foisf}
x(k+1)=x(k)+\varepsilon f(k,x(k))  
\ee
denote its Euler discretization, with~$\varepsilon>0$.

\begin{Lemma}\label{lem:edist}
Suppose that the trajectories of~\eqref{eq:foisf} evolve on a compact and convex set~$\Omega\subset\R^n$, and that
\be\label{eq:jacodis}
J(k,a) := I+ \varepsilon \frac{\partial}{\partial x}f(k,a)
\ee
 is tridiagonal, with positive entries on the super- and sub-diagonals
 for all~$k\in \mathbb{N}$ and all~$a\in \Omega$. Then  for 
any~$\varepsilon>0$ sufficiently small
Assumption~\ref{assume:jaco1} holds with~$h=n-1$. 
\end{Lemma}

Note that since~$J(k,a)$ is tridiagonal, it is   not~TP,
so the~TPDTS framework cannot be used to analyze this case. 

\begin{IEEEproof}
Pick~$k\in\mathbb{N}$ and~$a\in \Omega$. The assumptions on the Jacobian
 imply 
that~$J(k,a)$ is irreducible for all~$\varepsilon>0$, and
nonsingular for 
  every~$\varepsilon>0$ sufficiently small.
Also, $J(k,a)$ satisfies the dominance condition described in
Example~\ref{exa:osc} for 
  any~$\varepsilon>0$ sufficiently small, and is thus~TN.
 Furthermore, all these properties carry 
over to the matrix~$F$ defined in~\eqref{eq:deffmat1}. 
\end{IEEEproof}

The next example demonstrates 
Lemma~\ref{lem:edist} in a simple case.

\begin{Example}\label{exa:takac}
Consider the continuous-time  system:
\[
\dot x=
\begin{bmatrix} 0&1&0\\1&0&1\\0&1&0 \end{bmatrix} x.
\]
\end{Example}
Its Euler discretization  is~$x(k+1)=Ax(k)$ with
\[
A:=  \begin{bmatrix} 1&\varepsilon&0\\
\varepsilon&1&\varepsilon\\0&\varepsilon&1 \end{bmatrix},
\]
with~$\varepsilon>0$. 
The matrix~$A$ is irreducible, 
and it is nonsingular for
 any~$\varepsilon\not  =   1/\sqrt{2} $. 
Combining this with Example~\ref{exa:osc}
implies that~$A$ is oscillatory for
 any~$\varepsilon\in (0,1/2]$.  
The eigenvalues  of~$A$ are
\begin{align*}
\lambda_1&:=1+\varepsilon\sqrt{2} ,\\
\lambda_2&:=1,\\
\lambda_3&:=1-\varepsilon\sqrt{2} ,
\end{align*}
with corresponding eigenvectors
\begin{align*}
 v^1&:=\begin{bmatrix}1, \sqrt{2}, 1\end{bmatrix}',\\ 
v^2&:=\begin{bmatrix}-1, 0, 1\end{bmatrix}',\\
 v^3&:=\begin{bmatrix}1, -\sqrt{2}, 1\end{bmatrix}'.
\end{align*}
Pick~$x(0)\in\R^3$. Let~$c_i\in\R$ be such that
$
x(0)=\sum_{i=1}^3 c_i v^i.
$
Then
\[
x(k)=\sum_{i=1}^3 c_i \lambda_i^k v^i,
\]
and this implies that for
 any~$\varepsilon\in (0,1/2]$   
any solution~$x(k)$
 that remains in a compact set converges to either~$c_1 v^1$ or to the origin.

\begin{figure}
 \begin{center}
    \includegraphics[scale=0.6]{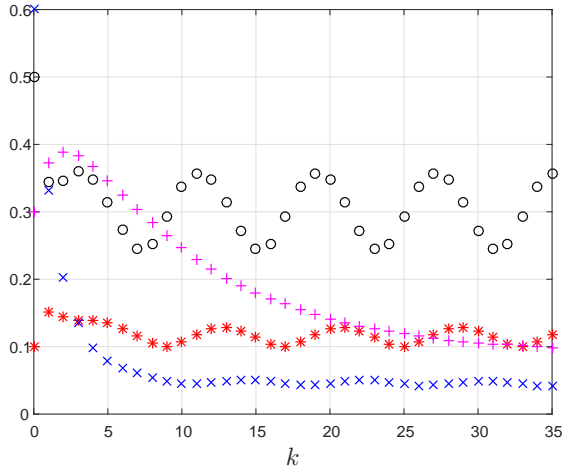}
	\caption{The trajectory~$x_1(k)$ (marked by o), $x_2(k)$~(*),
		$x_3(k)$~(x), and $x_4$~(+)   in Example~\ref{exa:phoso}. }	\label{fig:phos}
\end{center}
\end{figure}

 The next example describes an application of Lemma~\ref{lem:edist}
to a nonlinear 
 model from systems biology. 
\begin{Example}\label{exa:phoso}
Cells often sense and  respond to  various  stimuli  by
 modification of protein production. One   mechanism  for  this
 is \emph{phosphorelay}  (also called phosphotransfer),
 in which a phosphate group is transferred through a serial
 1D chain of proteins
  from an initial  histidine
kinase~(HK)    down to a
final response regulator~(RR).
The nonlinear compartmental system:
\begin{align}\label{eq:phos}
								\dot x_1&= (p_1-x_1)c-\eta_1x_1(p_2-x_2)-\xi_1 x_1,\nonumber\\
          			\dot x_2&=\eta_1x_1(p_2-x_2)-\eta_2x_2(p_3-x_3)-\xi_2 x_2,\nonumber\\
								&\vdots\nonumber\\
								\dot x_{n-1}&=\eta_{n-2} x_{n-2}(p_{n-1}-x_{n-1})-\eta_{n-1} x_{n-1}(p_n-x_n)\nonumber \\&-\xi_{n-1}x_{n-1},\nonumber\\
								\dot x_n&= \eta_{n-1} x_{n-1}(p_n-x_n) -\eta_n x_n,
\end{align}
has been suggested as a model for phosphorelay~\cite{phos_relays}.
Here~$c:[t_1,\infty)\to \R_+$ is the strength at time~$t$ of the
stimulus activating the~HK,
$x_i(t)$ is the concentration of
the phosphorylated form of the protein at the~$i$'th layer at time~$t$, the parameter~$p_i>0$ denotes the  total protein concentration at that layer,
and~$\eta_i,\xi_i>0$ are parameters that  describe reaction rates.
 Note that
$\eta_n x_n(t)$
is the  flow at time~$t$ 
of the phosphate group to an external  receptor molecule.

In  the particular case where~$p_i=1$ and~$\xi_i=0$
 for all~$i$ Eq.~\eqref{eq:phos}
  becomes the \emph{ribosome flow model}~(RFM)~\cite{reuveni}.
	This  is the dynamic mean-field approximation of a fundamental 
	model from non-equilibrium statistical physics called the \emph{totally asymmetric simple  exclusion process}~(TASEP) \cite{solvers_guide}.
  The RFM describes the unidirectional 
	flow along a chain of~$n$ sites. The state-variable~$x_i \in[0,1]$ describes
	the normalized occupancy at  site~$i$, where~$x_i=0$ [$x_i=1$]
means that site~$i$ is completely free [full], and~$\eta_i$ is the capacity of the link that connects site~$i$ to site~$i+1$. This has been used to model and analyze
mRNA translation (see, e.g.,~\cite{rfm_max,RFM_model_compete_J,rfm_nets,RFM_r_max_density}), where every site corresponds to a group of
codons on the mRNA strand,~$x_i(t)$
is the normalized occupancy   of ribosomes at site~$i$ at time~$t$,
$c(t)$ is the initiation rate at time~$t$, and $\eta_i$
 is the   elongation rate from site~$i$ to site~$i+1$.

Write~\eqref{eq:phos} as~$\dot x=f(x)$.
Then~$\frac{\partial}{\partial x}f(x)$
 is tridiagonal, with entries~$\eta_i x_i$  on the super-diagonal,
and~$\eta_i(p_{i+1}-x_{i+1})$, $i=1,\dots,n-1$, on the sub-diagonal. 

Consider the corresponding  discretized
 system~\eqref{eq:foisf}. It is not difficult to show that~$\Omega:=[0,p_1]\times\dots\times[0,p_n]$ is an invariant set of~\eqref{eq:foisf} for any~$\varepsilon>0$ sufficiently small. Furthermore, for any~$a\in\Omega$ we have that~$x(k,a)\in\Int(\Omega)$ for all~$k\geq 1$ and then 
the conditions in Lemma~\ref{lem:edist}
on~$J(k,a)$ defined in~\eqref{eq:jacodis}
hold.
Fig.~\ref{fig:phos} depicts the trajectories of the discretized system with~$n=4$, $\varepsilon=0.1$,
$\xi_i=3$,
$\eta_i=1$, $p_1=0.8$, $p_2=p_3=p_4=2$, 
initial condition~$x(0)=\begin{bmatrix}0.5& 0.1& 0.6& 0.3 \end{bmatrix} '$, 
and the periodic stimulus~$c(k)=3+\sin(k \pi/4)$.
      Note that this means that the map is~$T$-periodic with (minimal)
			period~$T= 8$. Combining Thm.~\ref{Thm:Oscillsys}
			and Lemma~\ref{lem:edist}, we conclude that any solution of 
			the discretized system converges
			to a periodic solution with period~$(n-1)T=24$. 
			It may be seen that the specific solution depicted in Fig.~\ref{fig:phos}
			converges to a periodic solution with period~$8$.
\end{Example}

In general, our approach is to
  find sufficient  conditions guaranteeing that the line integral
	of a matrix is oscillatory 
without actually calculating the integral. However, there is an  
important special case where
the integral can be computed explicitly.

\subsection{The case of strictly  monotone   scalar nonlinearities}
Let~$f_i:\R\to\R$,~$i=1,\dots,n$,
 be~$C^1$ functions such that 
\be\label{eq:yder}
f'_i(y):=\frac{d}{dy}f_i(y)>0 \text{ for all } i
\text{ and all } y\in\R.
\ee
 
 Consider the  time-varying
nonlinear system:
\be\label{eq:dofty}
x(k+1)=C(k) \begin{bmatrix} f_1(x_1(k))  
\\f_2(x_2(k))\\\vdots\\f_n (x_n(k))\end{bmatrix},
\ee
	with~$C:\N\to\R^{n\times n }$.	
	
\begin{Theorem} \label{thm:repnon}
Suppose that the trajectories of~\eqref{eq:dofty}
evolve on a compact and 
convex state-space~$\Omega$,
and that~$C(k)$ is~$T$-periodic.  
If~$z(k+1)=C(k)z(k)$ is an~ODTS of order~$h$
  then 
	 every solution of~\eqref{eq:dofty} emanating from~$\Omega$ 
converges to an~$hT$-periodic solution of~\eqref{eq:dofty}.
\end{Theorem} 

\begin{IEEEproof}
  The Jacobian of~\eqref{eq:dofty}  is
$
J(k,x)= C (k)  \diag(f'_1(x_1)),\dots,f'_n(x_n))  
$.
Substituting this in~\eqref{eq:deffmat} and integrating
 yields
\begin{align}\label{eq:jacoiop}
                F(k,a,b)
                &= C (k) \diag(g_1(a_1,b_1) ,\dots,g_n(a_n,b_n) ),
\end{align}
   with
\begin{align*}
 g_i(a_i,b_i) 
&:= \begin{cases}   \frac{f_i(a_i)-f_i(b_i)}{a_i -b_i} & \text{ if }  
a_i\not= b_i,\\
       f'_i(b_i)  & \text{ if } a_i = b_i.                                \end{cases}
\end{align*}
Note that~\eqref{eq:yder} and the fact that~$\Omega$ is compact imply that 
there exists~$\delta>0$ such that~$g_i(a_i,b_i)\geq \delta$ for all~$a,b\in \Omega$ and all~$i$.

Pick~$1\leq  r \leq n$, and indexes~$1\leq i_1 < \dots< i_r\leq n$ and~$1\leq j_1 <  
\dots < j_r\leq n$.
Let~$F(\alpha|\beta)$ denote the minor of~$F=F(k,a,b)$ indexed by rows~$i_1,\dots,i_r$
and columns~$j_1,\dots,j_r$. Then applying the Cauchy-Binet formula  (see, e.g \cite{total_book})
to~\eqref{eq:jacoiop}
yields
\be\label{eq:fmino}
F(\alpha|\beta)=C(\alpha|\beta) g_{j_1} g_{j_2} \dots g_{j_r}.
\ee
Since all the~$g_i$'s are positive, this means that
the total positivity properties of~$C$
are copied to~$F$. 
  Applying Thm.~\ref{Thm:Oscillsys} completes the proof.
\end{IEEEproof}

\begin{figure}
 \begin{center}
    \includegraphics[scale=0.6]{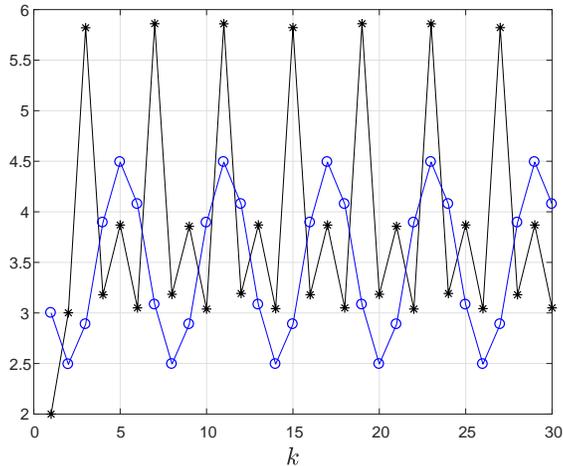}
	\caption{The trajectory~$x_1(k)$ (marked by *), $x_2(k)$ (o)  for the initial condition~$x_1 (0) =2,\; x_2 (0) = 3$ in Example~\ref{exa:capc}. }	\label{fig:entrainment}
\end{center}
\end{figure}

The next example demonstrates Thm.~\ref{thm:repnon}.
\begin{Example}\label{exa:capc}
Consider the system:
\begin{align}\label{eq:bfty}
\begin{bmatrix}
x_1 (k+1) \\
x_2 (k+1)
\end{bmatrix}
=
C(k)\begin{bmatrix}
\tanh(x_1(k)) \\ \tanh(x_2(k)) 
\end{bmatrix},
\end{align}
where
\begin{align*}
c_{11}(k) &= 2+\cos(k\pi+0.5),\\ 
c_{12}(k) &=  2-\sin(\frac{k\pi}{2}+1.5) ,\\
c_{21}(k) &\equiv  1/2 ,\\
c_{22}(k)& = 3+\cos(\frac{k\pi}{3} +2).
\end{align*}
 Note that~$C(k)$ is~TP for all~$k\in\N$,  and that
\be\label{eq:cbo}
1\leq |c_{ij}(k)|\leq 4 \text{ for all } i,j,k,
\ee
and  that  the map in~\eqref{eq:bfty}
  is periodic with (minimal) period~$T=12$. 

We claim that (for example)  the square
\[
\Omega:=[1,8]\times [1,8]
\]
is an invariant set for the dynamics. To show this, suppose that~$x(k)\in \Omega$. Then~$x_1(k),x_2(k)\geq 1$, 
so~$x(k+1)\geq\begin{bmatrix}1&1\\1/2&2 \end{bmatrix}  \begin{bmatrix}\tanh(1)\\ \tanh(1)\end{bmatrix}
\geq \begin{bmatrix}1\\1\end{bmatrix}$. 
Also, since~$\tanh(z)\leq 1$ for all~$z $, and~$c_{ij}(k)\leq 4$,  
we have~$x_1(k+1),x_2(k+1)\leq 8$, so~$x(k+1)\in \Omega$.  

It is clear that this system satisfies the conditions in Thm.~\ref{thm:repnon}, 
with~$h=1$, and thus the system entrains. 
The trajectory of the system
for~$x_1 (0) = 2$ and  $x_2 (0) = 3$ 
is depicted 
  in Fig.~\ref{fig:entrainment}.
It may be seen that~$x_1 (k), x_2 (k)$ indeed 
converges to a
 $T$-periodic solution with~$T=12$.  
\end{Example}

From here on we consider 
 the following general  problem.

\begin{Problem}\label{prob:wje}
Consider a  measurable and essentially bounded
 matrix function~$A:[0,1]\to\R^{n\times m}$. 
When is
\[
\bar A :=\int_0^1   A(t)  \diff t
\]
  an oscillatory  matrix?
\end{Problem}
 
 Since our motivation is the analysis
of  line integrals of   Jacobians of dynamical systems, 
we assume throughout that~$m=n$. Some of the conditions given
below actually guarantee that~$\bar A$
  is~TP (and thus, in particular, oscillatory with exponent one).
\subsection{Sufficient condition based on  the checkerboard partial order}
For  $A,B \in \mathbb{R}^{n\times n}$ we write $A\leq B$ [$A\ll B$]
 if $a_{ij}\leq b_{ij}$  [$ a_{ij}<b_{ij}$]
 for all $i,j $. 

\begin{definition}
 The checkerboard partial order on $\mathbb{R}^{n\times n}$
is defined by 
\[
A \leq^{\dagger}B \iff
 D_\pm A D_\pm \leq D_\pm B D_\pm .\]
\end{definition}


In other words,  $A\leq^{\dagger}B$ iff 
\be\label{eq:aldef}
(-1)^{i+j}a_{ij}\leq (-1)^{i+j}b_{ij} \text{ for all } i,j \in\{1,\dots,n\}.
\ee 
Note that~\eqref{eq:aldef} implies that 
the  matrix interval~$\left\{C \in \mathbb{R}^{n \times n}: A\leq^{\dagger}C\leq^{\dagger}B\right\}$ is compact. For more on such matrix
intervals, see~\cite{GARLOFF2003103} and the references therein.

It is well-known~\cite{total_book} that 
if~$A,B$ are~TP and~$A\leq^{\dagger}C\leq^{\dagger}B$ then $C$ is~TP.

\begin{Theorem}\label{Thm:TPint1}
Let $A:[0,1]\to\R^{n\times n}$ be a Riemann
  integrable  
 matrix function. 
If there exist~$\delta>0$ and~TN matrices~$G$ and~$H$ such that
\begin{eqnarray}\label{eq:TPline}
\delta+(-1)^{i+j} g_{ij} 
 \leq (-1)^{i+j} a_{ij}(t)\leq - \delta   + (-1)^{i+j} h_{ij}   
\end{eqnarray}
  for all $i, j$  and all~$t \in   [0,1]$ 
then~$\bar A$ is TP.
\end{Theorem}

\begin{IEEEproof} 
Recall  that the set of~$n\times n$
TP matrices  is dense in the set of  $n \times n$ TN
 matrices~\cite{Whitney1952}.
Combining  this with~\eqref{eq:TPline} implies that there exist~TP matrices~$P$ and~$Q$ such that
\be\label{eq:bnfpo}
P\leq^{\dagger}A(t) \leq^{\dagger}Q \text{ for all } t\in[0,1] .
\ee
We claim  that this implies that every  minor of~$\bar A$ is positive. 
We will show that~$\det \bar A >0$. The proof for any other minor is very similar. 
Fix~$k\in \{1,2\dots,\}$ and consider the partition of~$[0,1]$  
 defined by
\[
t_0:=0, \; t_1:=1/k, \;  t_2:=2/k,\dots,t_k:=1.
\]
Consider the Riemann   sum 
  $B:= \sum_{\ell=0}^{k-1} (t_{\ell+1}-t_\ell)A(t_\ell)$. Then for any~$i,j\in\{1,\dots,n\}$ we have
	\[
	(-1)^{i+j} b_{ij} = \sum_{\ell=0}^{k-1} (-1)^{i+j} (t_{\ell+1}-t_\ell)a_{ij}(t_\ell),
	\]
and combining this with~\eqref{eq:bnfpo} gives 
\begin{align*}
	 \sum_{\ell=0}^{k-1} (-1)^{i+j} (t_{\ell+1}-t_\ell)p_{ij}&\leq (-1)^{i+j} b_{ij}\\
	&\leq \sum_{\ell=0}^{k-1} (-1)^{i+j} (t_{\ell+1}-t_\ell)q_{ij}.
\end{align*}
Since~$\sum_{\ell=0}^{k-1}  (t_{\ell+1}-t_\ell) =t_k-t_0=1$, 
we conclude that 
\[ 
P\leq^{\dagger}B \leq^{\dagger}Q .
\]
By compactness of the set~$\{ C\in\R^{n\times n}: P\leq^{\dagger} C \leq^{\dagger}Q \}$ and the fact that any~$C$ in this set is~TP,  there exists~$\alpha>0$ such that~$\det B \geq \alpha$. Taking~$k \to\infty$ and using 
the continuity of the determinant, we conclude  that~$\det \bar A \geq\alpha >0$.
\end{IEEEproof}

Suppose  that   every 
entry~$a_{ij}(t)$ of~$A(t)$
attains a maximum value~$\bar a_{ij}$
  and a   minimum~$\underline a_{ij}$ over~$[0,1]$.
	Define~$P, Q$ by
\[
p_{ij}:=\begin{cases}   
\underline a_{ij}, & \text{ if } i+j \text{ is even}, \\
\bar a_{ij}, & \text{ if } i+j \text{ is odd},  
\end{cases} 
\]
and
\[
q_{ij}:=\begin{cases}   
\bar a_{ij}, & \text{ if } i+j \text{ is even}, \\
\underline  a_{ij}, & \text{ if } i+j \text{ is odd},  
\end{cases} 
\]
 Then~\eqref{eq:bnfpo} holds, so the required condition is that~$P$ 
and~$Q$ are~TP.

The next result
 describes
 an application of Thm.~\ref{Thm:TPint1} to 
a dynamical system.
\begin{Corollary}
Consider the nonlinear system:
\be\label{eq:pert}
x(k+1)=Ax(k)+\varepsilon g(x(k)),
\ee
where~$g$ is~$C^1$ and~$\varepsilon>0$ is small.
Suppose that~$A$ is TP,
and that 
the trajectories of~\eqref{eq:pert}
evolve on a compact and
 convex set~$\Omega\subset\R^n$. Define~$B\in \R^{n\times n}$ by
\[
b_{ij}:=\max_{x\in \Omega} \left |\frac{ \partial g_{i } (x) }{\partial x_j  } \right |,
\]
and define
 matrix functions~$P,Q:\R\to\R^{n\times n}$ by
\be\label{eq:defpq}
P(v):=A-v D_\pm B D_\pm,\quad Q(v):=A + v D_\pm B D_\pm.
\ee
Then there exists~$w>0$ such that for all~$v\in[0,w)$
 and all~$x\in \Omega$ 
\[
      P(v) \leq^\dagger   J(x) \leq^\dagger  Q(v),
\]
and~$P(v),Q(v)$ are TP, and for any~$\varepsilon \in[0,w)$ 
every solution of~\eqref{eq:pert}  emanating from~$\Omega$ 
converges to an equilibrium point. 
\end{Corollary}

\begin{IEEEproof}
It follows from~\eqref{eq:defpq} that~$P(0)=Q(0)=A$,
 and~$P(v)\leq^\dagger A \leq^\dagger Q(v)$ for all~$v\geq 0$. 
By continuity of the minors, there exists~$w>0$ such that
\[
P(v),Q(v) \text{ are TP for all } v\in [0,w). 
\]

 The Jacobian of~\eqref{eq:pert} is
$
J(x )=A+\varepsilon \frac{\partial}{\partial x} g(x ), 
$
so for any~$s,r\in\{1,\dots,n\}$ and any~$x\in\Omega$ we have 
\begin{align*}
|J_{sr}(x)|&=| a_{sr}+\varepsilon  \frac{\partial}{\partial x_r} g_s(x)| \\
&\leq a_{sr}+\varepsilon b_{sr}.
\end{align*}
It is straightforward to verify that this implies that for any~$v\geq 0$
and any~$\varepsilon \in [0,v]$ we  have  
\[
(-1)^{s+r} p_{sr}(v) \leq (-1)^{s+r}  J_{sr}(x) \leq (-1)^{s+r} q_{sr}(v),
\]
that is,
\[
P(v)\leq^\dagger J(x) \leq^\dagger Q(v).
\]
Now  fix~$\varepsilon \in[0,w)$. Pick~$v\in[\varepsilon,w)$.
Then for these values all the conditions in 
Thm.~\ref{Thm:TPint1} hold, so 
the matrix~$F(a,b)$ in~\eqref{eq:defab} is~TP for all~$a,b\in\Omega$,
and this completes the proof.
\end{IEEEproof}

\begin{figure}[t]
 \begin{center}
  \includegraphics[scale=0.5]{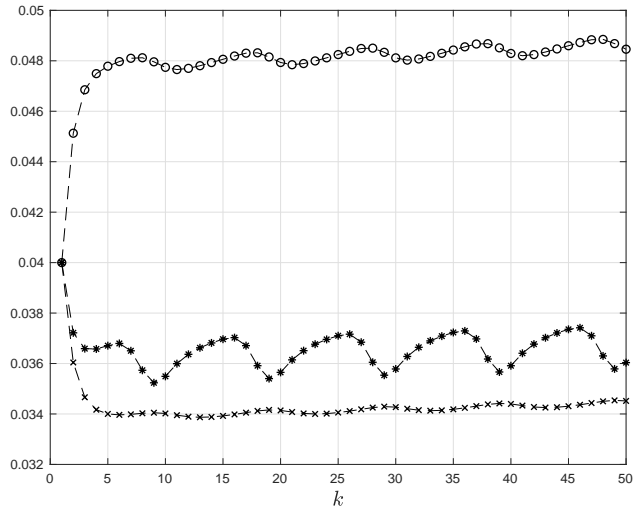}
\caption{State-variables $x_1(k)$ (marked with~'*'), $x_2(k)$ ('o'),
and~$x_3(k)$ ('x')  as a function of~$k$
for the system in Example~\ref{exa:pertu}. }
\label{fig:secent}
\end{center}
\end{figure}

\begin{Example}\label{exa:pertu}
Consider~\eqref{eq:pert} with~$n=3$, 
\be\label{eq:abrew}
A=0.65 \begin{bmatrix}
1& \exp(-1) &\exp(-4)\\
\exp(-1) &1&\exp(-1) \\
\exp(-4) &\exp(-1) & 1
\end{bmatrix},
\ee
and~$g(k,x(k))=\begin{bmatrix} \tanh ( (50+50\sin(k\pi/5)) x_3(k) )  &0&0  \end{bmatrix}
'$. 
This model 
may represent a cooperative 
linear chain where the effect of~$x_i(k)$ on~$x_j(k+1)$
decays exponentially with the ``distance''~$(i-j)^2$
between~$x_i$ and~$x_j$. 
It is well-known that~$A$ in~\eqref{eq:abrew}
 is~TP (see~\cite[Ch.~II]{gk_book}).  The nonlinear term 
represents a time-varying and~$T$-periodic, with~$T=10$,
positive feedback from~$x_3$ to~$x_1$.
 
It is clear that we can take the ``bounding matrix''~$B\in \R^{3\times 3}$
as the matrix with~$b_{13}=1$, and zero in all other entries. It is not difficult to verify that for this~$B$ we have 
that~$P(v),Q(v)$ defined in~\eqref{eq:defpq}
 are~TP
for all~$v\in[0,w)$, with~$w:=0.65\exp(-4)$. Fig.~\ref{fig:secent}
depicts the solution of the system with~$\varepsilon=0.0118<w$
and initial condition~$x(0)=(2/50) \begin{bmatrix} 1  & 1    &1 \end{bmatrix}'$. It may be seen that every~$x_i(k)$ converges to a periodic solution with period~$T=10$.  
\end{Example}

\subsection{Integrating TP Hankel Matrices}

Recall that~$A\in\R^{n\times n}$ is called a \emph{Hankel matrix} if
for any~$i,j,p,q$ with~$i+j=p+q$ we have~$a_{ij}=a_{pq}$.  For example, for~$n=3$
a Hankel matrix has the form
\[
A=\begin{bmatrix} 
a_{11}&a_{12}&a_{13} \\
a_{12}&a_{13}&a_{23} \\
a_{13}&a_{23}&a_{33} 
\end{bmatrix}.
\]
Note that a Hankel matrix is in particular 
 symmetric. Our main result in this subsection is
that the integral of a time-varying TP Hankel matrix is~TP.
\begin{Theorem}\label{Thm:HankelIntegral}
Let $A:[0,1]\to\R^{n\times n}$ be a measurable matrix function such that~$A(t)\in L^{\infty}([0,1])$. Suppose that~$A(t)$ is a~TP Hankel matrix for almost every $t\in [0,1]$. Then~$\bar A$ is TP.
\end{Theorem}

\begin{remark}\label{rem:n22}
Note that for~$n=2$ this implies that if~$A:[0,1]\to\R^{2\times 2}$ is  a continuous matrix function with~$A(t)$   symmetric and~TP
   for all~$t\in [0,1]$ then~$\bar A$ is~TP (compare with Example~\ref{exa:nirt}).
\end{remark}

To prove Thm.~\ref{Thm:HankelIntegral}
 we recall several definitions and results. 
A set of indices~$I\subseteq \left\{1,\cdots ,n\right\}$ is called an \emph{interval}
if it has the form~$I=\{p,p+1,p+2,\dots, q\}$.
A square sub-matrix  of a matrix~$B\in\R^{n\times n}$
 with row indices~$I \subseteq \left\{1,\cdots ,n\right\}$ and column indices~$J \subseteq \left\{1,\cdots ,n\right\}$ is called a \emph{contiguous sub-matrix}
 if both~$I$ and~$J$ are   intervals. 

It is well-known and straightforward to show
 that the following three conditions are equivalent:
(1)~$B \in \R^{n\times n}$  is a Hankel matrix;
(2)~every contiguous sub-matrix of $B$ is a Hankel matrix;
(3)~every contiguous sub-matrix of $B$ is symmetric.
 
 We can now prove Thm.~\ref{Thm:HankelIntegral}.

\begin{IEEEproof} 
We start by showing that $\det \bar A >0$. First, note that the function $t \mapsto \det( A(t))$ is measurable (as it is a polynomial in the entries~$a_{i,j}(t)$, $i,j \in \left\{1,\dots,n\right\}$) and essentially bounded. Therefore, it is Lebesgue integrable. For~$N\in \{1,2,\dots,\}$, let
\begin{equation*}
B_N:= \left\{t\in[0,1] : \det( A(t)) \geq N^{-1}\right\}.
\end{equation*}
Since $A(t)$ is~TP for almost every~$t\in [0,1]$ and \[
B_1\subseteq B_2 \subseteq \dots,
\]
   the monotone convergence theorem
 (see e.g.~\cite{bogachev2007measure}) yields 
\begin{equation*}
\lim_{N \rightarrow \infty}\mu(B_N) =1,
\end{equation*}
where $\mu$ is the Lebesgue measure on~$[0,1]$. Therefore, there exists~$N_0 \in \mathbb{N}$ such that $\mu(B_{N_0})>\frac{1}{2}$. 
Markov's inequality (see e.g.~\cite{bogachev2007measure}) yields 
\begin{align}\label{eq:Markov}
\int_{0}^{1} \left ( \det{A(t)} \right) ^{\frac{1}{n}}\diff \mu(t) &\geq N_0^{-\frac{1}{n}}\mu(B_{N_0})\nonumber\\
&>  N_0^{-\frac{1}{n}}/2.
\end{align}

Since~$A(t)$ is  Hankel and~TP 
for almost all~$t\in[0,1]$, it is symmetric with positive principal minors, so~$A(t)$   is positive-definite for almost all~$t\in[0,1]$.
Minkowki's determinant inequality (see e.g. \cite[p. 115]{marcus1992survey}) states
 that~$B \mapsto (\det B)^{\frac{1}{n}}$ is a concave function over
 the space of semi-positive definite matrices of order~$n$. Thus, by using \eqref{eq:Markov} and Jensen's inequality (see e.g.~\cite{bogachev2007measure}) we obtain
\begin{align*}
(\det \bar A )^{\frac{1}{n}}&=\left( \det \int_0^1   A(t)  \diff \mu(t)\right)^{\frac{1}{n}} \\
&\geq \int_0^1   (\det A(t))^{\frac{1}{n}}  \diff \mu(t) \\
& \geq  N_0^{-\frac{1}{n}}/2,
\end{align*}
so~$\det \bar A >0$. 

Recall that every contiguous sub-matrix of $A(t)$ is also a~TP and
Hankel matrix for almost all~$t\in [0,1]$, so the same argument shows that every 
contiguous minor of $\bar A$ is positive.
It is well-known~\cite[Chapter~3]{total_book}
 that if all the contiguous minors of a  
 matrix are positive then the matrix is~TP, so we conclude that~$\bar A$ is~TP.
\end{IEEEproof}

The next example demonstrates an application of 
Remark~\ref{rem:n22}
to a dynamical system.
\begin{Example}
Consider  the   nonlinear system:
\begin{align}\label{eq:frvce}
x_1(k+1)&=h_1(x_1(k))+g(x_1(k),x_2(k)), \nonumber \\
x_2(k+1)&=h_2(x_2(k))+g(x_1(k),x_2(k)),
\end{align}
with~$h_1,h_2,g\in C^1$,
whose trajectories evolve on a compact and convex state-space~$\Omega\subset\R^2$. 
Suppose that~$\frac{\partial}{\partial x_1}g(x_1,x_2)=\frac{\partial}{\partial x_2}g(x_1,x_2)$
for all~$x_1,x_2\in \Omega$ (e.g.~$g(x_1,x_2)=\tanh(x_1+x_2)$). Note that this implies that the Jacobian
\[
J (x )= \begin{bmatrix}
h_1'(x_1) +\frac{\partial}{\partial x_1}g(x_1,x_2) & \frac{\partial}{\partial x_2}g(x_1,x_2)\\
 \frac{\partial}{\partial x_1}g(x_1,x_2) &h_2'(x_2)+ \frac{\partial}{\partial x_2}g(x_1,x_2)\\
\end{bmatrix}
\] 
is symmetric. If~$J(x_1,x_2)$ is~TP for all~$(x_1,x_2)\in\Omega$
then combining Corollary~\ref{coro:tinv} and
Remark~\ref{rem:n22}
  implies that any solution of~\eqref{eq:frvce}
emanating from~$\Omega$ converges to an equilibrium point.
\end{Example}

\section{Conclusion}\label{sec:conc}

We introduced a new class of positive  discrete-time LTV
  systems called~ODTSs of order~$p$. Discrete-time
nonlinear systems, whose variational system is 
an~ODTS of order~$p$, have a well-ordered behavior.
More precisely, if the map defining the dynamical 
system is~$T$-periodic then every solution either
leaves any compact set or converges to a~$pT$-periodic solution, i.e. a  subharmonic
  solution.
 This is important because, as 
noted by 
Smith~\cite{smith_planar}, 
``\dots in the
class of all discrete dynamical systems, we do not know so many
special classes which have relatively simple dynamics.''

The  ODTS framework  
requires establishing  that certain
 line integrals of the Jacobian of the time-varying
nonlinear system   are  oscillatory matrices. 
This is non-trivial, as the sum of two oscillatory matrices is not necessarily oscillatory,  and this naturally extends to integrals.
We derived several sufficient conditions
 guaranteeing that the line integral of a matrix is~oscillatory (or~TP).

Topics for further research include the following. 
First, extending the oscillatory framework to other dynamical models e.g. systems with time-delays
 or  discretized~PDEs. Second, cooperative
discrete-time systems arise frequently as the Poincar\'e 
maps of continuous-time systems.  It may be of interest   to explore the implications of   oscillatory 
Poincar\'e  maps. Third, it may be of interest to  generalize
 the~ODTS framework to discrete-time systems with control inputs, as was done for continuous-time monotone
 systems in~\cite{mcs_angeli_2003}.

\section*{Appendix}
 \begin{IEEEproof}[Proof of Lemma~\ref{lem:nrety}]
Let~$z:=D_\pm (  x-  y)$. Then~\eqref{eq:pozie} implies
 that~$v:= D_\pm Pz\ll 0$. Thus,
  the vector~$Pz= D_\pm v $ is alternating, with
\[	(Pz)_1=v_1<0.
\]
	 Applying  the 
   VDP~\eqref{eq:vdptn} yields
	\[
						n-1=	s^-( P z ) \leq s^-(   z ) .
	\]
 Thus,~$s^-(  z ) =n-1$,  i.e.~$  z$ is alternating. 
Recall that if a matrix~$H$ is~TN, and non-singular 
and~$s^-(Hq)=s^-(q)$ for some~$q\in\R^n\setminus\{0\}$
then the first non-zero entry in~$Hq$ and 
the  first non-zero entry in~$q$
have the same sign~\cite[p.~254]{gk_book}. 
Since~$s^-(P  z)=s^-(  z)=n-1$, and~$(Pz)_1<0$, 
the first non-zero entry of~$  z$  is negative. 
Since~$z$ is alternating 
 this implies that~$D_\pm   z\ll 0$,
 and this completes the proof.
\end{IEEEproof}


\end{document}